\documentclass{amsart}
\usepackage{latexsym,amssymb,amsfonts,amsmath, graphicx, color,enumerate}
\usepackage{mathbbol}
\usepackage{dsfont}
\usepackage[latin1]{inputenc}
\usepackage{hyperref}

 \newtheorem{thm}{Theorem}[section]
 \newtheorem{cor}[thm]{Corollary}
 \newtheorem{lem}[thm]{Lemma}
 \newtheorem{prop}[thm]{Proposition}
 \theoremstyle{definition}
 
 \theoremstyle{remark}
 \newtheorem{rem}[thm]{Remark}
 
 \numberwithin{equation}{section}

\def\:{\thinspace:\thinspace}
\def\ea{\mathbf a}
\def\eu{\mathbf u}
\def\ev{\mathbf v}
\def\ez{\mathbf z}
\newtheorem{exa}[thm]{Example}

\numberwithin{equation}{section}
\numberwithin{thm}{section}
\def\:{\thinspace:\thinspace}
\begin{document}

\title[A variational approach to strongly damped wave equations]{A variational approach to\\ strongly damped wave equations}

\author{Delio Mugnolo}

\address{Institut f\"ur Analysis\\
Helmholtzstra{\ss}e 18\\
Universit\"at Ulm\\
D-89081 Ulm, Germany}
\email{delio.mugnolo@uni-ulm.de}

\thanks{This is the slightly enhanced version of an article published in the proceedings of the conference in memory of Gunther Lumer held in Mons and Valenciennes in 2006. The published version of this article lacks Section 4.}

\begin{abstract}
We discuss a Hilbert space method that allows to prove analytical well-posedness of a class of linear strongly damped wave equations.
The main technical tool is a perturbation lemma for sesquilinear forms, which seems to be new. In most common linear cases we can furthermore apply a recent result due to Crouzeix--Haase, thus extending several known results and obtaining optimal analyticity angle.
\end{abstract}

\maketitle

\section{Introduction}\label{intro}

Of concern of this note are complete second order abstract Cauchy problems of the form
\begin{equation}\label{damped}
\left\{\begin{array}{l}
\ddot{u}(t)+Au(t)+B\dot{u}(t)=0,\qquad t\geq 0,\\
u(0)=u_{10},\qquad \dot{u}(0)=u_{20},
\end{array}\right.
\end{equation}
where the elastic operator $A$ is in the literature usually assumed to be a self-adjoint, strictly positive definite operator on a Hilbert space $H$. It is known that such elastic systems exhibit good properties whenever $B$ is a multiplication operator: e.g., they are forward as well as backward solvable, they admit energy decay estimates if $B$ is dissipative, or else blow-up estimates if $B$ is accretive, see e.g.~\cite{KL78,KL78b,LG97} and references therein.

It is interesting to note that, in particular, the standard model of an electrical transmission line by means of the telegraph equation fits this framework, the case of $B$ negative multiplication operator corresponding to \emph{viscous} damping.

\bigskip	
In~\cite{CR82}, Chen--Russell proposed a family of different, strongly (or \emph{structural}) damping effects: theoretical arguments and empirical studies motivated them to consider damping operators that are unbounded on $H$, cf. references in~\cite{CR82}. For the sake of simplicity, they mostly investigated the special cases of $B=A$ and $B=2\rho A^\frac{1}{2}$. However, they also pointed out that the crucial property is the so-called \emph{frequency response} estimate
\begin{equation*}\label{sector}
\| \lambda R(i \lambda,{\mathbf A})\|\leq {M},\qquad \lambda\in\mathbb R,
\end{equation*}
satisfied by the resolvent operator of $\mathbf A$, where
\begin{equation}
 {\mathbf A}:=\begin{pmatrix}0 & -I\\ A & B\end{pmatrix}
 \end{equation}
is the reduction matrix associated with~\eqref{damped}. %In fact the property~\eqref{sector}, shared by all sectorial operators $\mathbf A$, ensures that it is possible to develop a spectral theory which allows for an effective modelling by finite-dimensional systems. In other words, 
Thus, following Chen--Russell the issue becomes to find conditions on $A,B$ ensuring that $\mathbf A$ (or rather its closure) generates an analytic semigroup in the candidate phase space ${\mathbf H}:=D(A^\frac{1}{2})\times H$. %Observe that while general damped wave equations do not admit a unique phase space (in sharp contrast to undamped ones), in concrete cases

Ever since, several authors including Dautray--Lions, Chen--Triggiani, Xiao--Liang, and Chill--Srivastava have further investigated these kind of parabolic systems, significantly extending the results of Chen--Russell. Chen--Triggiani still imposed the assumption that the damping effect is at most as strong as the elastic one, i.e., that
\begin{equation}\label{standing}
B=\rho A^\alpha,\qquad \hbox{for }\alpha\in [0,1]\hbox{ and }\rho\in(0,\infty),
\end{equation}	
and then showed, by methods based on spectral analysis,  that the semigroup generated by the closure (of a suitable part) of $-\mathbf A$ is analytic if and only if $\alpha\in[\frac{1}{2},1]$, cf.~\cite[Thm.~1.1]{CT89}.  Successively, Xiao--Liang have proved similar results in the slightly more general case where $B=f(A)$ for a suitable class of functions $f$%. Roughly speaking, they proved that if $f(x)=o(x^\alpha)$, $\alpha\in(\frac{1}{2},1]$, then the semigroup generated by the closure (of a suitable part) of $\mathbf A$ is analytic of angle $\frac{\pi}{2}$
, cf.~\cite[Thm.~6.4.2]{XL98}. Similar, less sharp results have also been obtained in~\cite[\S~6.3]{EN00} by a technique based on the theory of operator matrices. We observe that strongly damped wave equations are also of interest in the framework of control theory, see e.g.~\cite{Bu93,LPT97}, and references therein. Energy decay estimates have also been extensively investigated, see e.g.~\cite{Hu97,BE04}.

% More recently, Chill--Srivastava have investigated similar systems by completely different methods based on the theory of $L^p$-maximal regularity, as they were mainly interested in properties of inhomogeneous problems of the form
%In fact, having established $L^p$-maximal regularity for~\eqref{dampedinh}, it is almost standard to obtain well-posedness and regularity results for more general classes of semilinear or even fully nonlinear problems. 
More recently, Chill--Srivastava have discussed  $L^p$-maximal regularity properties for the solution to
\begin{equation}\label{dampedinh}
\left\{\begin{array}{l}
\ddot{u}(t)+Au(t)+B\dot{u}(t)=f(t),\qquad t\in [0,T],\\
u(0)=0,\qquad \dot{u}(0)=0.
\end{array}\right.
\end{equation}
While they are not directly interested in the analyticity of the semigroup generated by $-\mathbf A$, their results in some sense extend those of~\cite{CT89,XL98}: if~\eqref{standing} holds and $A$ is a sectorial operator on an $L^q$-space, $q\in(1,\infty)$, and under further technical assumptions, it turns out that~\eqref{dampedinh} has maximal $L^p$-regularity if $\alpha\in (\frac{1}{2},1]$, cf.~\cite[Thm.~4.1]{CS05}. Observe that, in particular, if~\eqref{dampedinh} has $L^p$-maximal regularity and hence the evolution of the systems shows a smoothing effect, it is in general not clear how to formulate the second order problem in terms of a first order problem on a product space, i.e., it is not clear how to determine a natural phase space of the problem. Still if~\eqref{dampedinh} has $L^p$-maximal regularity, then a suitable phase space can actually be found in a natural way, cf.~\cite[Cor.~2.2]{CS05}.

The case of $D(B)\subset D(A)$ has been treated less frequently, see e.g.~\cite{Sa83,Ne86}; moreover, most authors have not discussed analyticity properties. In~\cite[Thm.~6.2]{Mu06b}, we have showed that if $B$ generates a cosine operator function with phase space $V\times H$, and if $A$ is bounded from $V$ to $H$, then~\eqref{damped} is governed by an analytic semigroup of angle $\frac{\pi}{2}$. In many relevant cases this amounts to saying that $A=\rho B^\alpha$, for $\alpha\leq\frac{1}{2}$ and $\rho\in\mathbb C$.

\bigskip
Aim of this paper is to discuss~\eqref{damped} under assumptions on $B$ that complement, or perhaps interpolate, those of the above mentioned papers. In fact, %while all the aforementioned authors assumed $B$ to be \emph{at most} as unbounded as $A$, 
we will assume $B$ to be \emph{at least} as unbounded as $A$. The quoted results suggest that $\alpha=\frac{1}{2}$ is a critical exponent, whenever~\eqref{standing} holds. In fact, we will show that the exponent $\alpha=1$ is critical, too. More precisely if $\alpha=1$, %or more generally if $D(A)=D(B)$, 
then the leading term in~\eqref{damped} is not $A$ anymore, but $B$. In fact, we show that~\eqref{damped} is governed by an analytic semigroup under quite weak boundedness assumptions on $A$, whenever $B$ is associated with a closed, $H$-elliptic form. In particular, we show that no closedness or spectral conditions on $A$ are necessary. Our method is based on the introduction of a suitable weak formulation 
%$$(\ddot{u}\mid \dot{u})+(Au\mid \dot{u})+(B\dot{u}\mid \dot{u})=0$$
of~\eqref{damped}, and then on the application of the theory of sesquilinear forms on complex Hilbert spaces. We refer to~\cite{Ou05,Ar04} for comprehensive treatments of this mature theory that goes back to Kato and Lions, and to~\cite{DL88} for a similar, slightly less general approach to damped wave equations due to Dautray--Lions.

In Section~\ref{first} we introduce our general framework and show a first  well-posedness result for~\eqref{damped}. To this aim we prove a perturbation lemma for sesquilinear forms that may be of independent interest.  We also obtain a first estimate on the angle of analyticity. In Section~\ref{interpo} we impose slightly stronger conditions and, by means of a recent result due to Crouzeix--Haase, we find sufficient condition in order  that the semigroup is analytic of angle $\frac{\pi}{2}$: this includes the relevant case of self-adjoint damping operator $B$. Some applications to semilinear problems are also considered. %	: this includes the relevant case of $B$ self-adjoint. %This in turn allows to discuss the well-posedness of~\eqref{dampedinh}, by means of the results presented in~\cite{Lu95}.
Finally, in Section~\ref{dynam} we briefly discuss how our theory can be adapted in order to discuss damped wave equations with a certain class of dynamic boundary conditions.

\bigskip
It is fair to add that a variational approach to linear damped wave equations has also been pursued in~\cite[\S~XVIII.5.1]{DL88}, see also~\cite[\S~XVIII.6]{DL88}. Indeed, Dautray--Lions' methods are quite similar to those presented in Section~2 below, and they also consider the neutral equation $C\ddot{u}(t)+Au(t)+B\dot{u}(t)=0$, $t\geq 0$, even in the nonautonomous case, where $D(B)\subset D(A)$. Though, the assumptions in~\cite[\S~XVIII.5.1]{DL88} are restricted to the case of $A,B$ differential operators whose principal part is self-adjoint and (in the case of $B$) also strictly positive definite, and no angle of analyticity is proved there. However, their main result~\cite[Thm.~XVIII.1]{DL88} is admittedly very close to Corollary~\ref{wellp} below.

\section{First well-posedness results}\label{first}

Let $V,H$ be separable complex Hilbert spaces such that $V$ is continuously and densely imbedded in $H$. Let $a:V\times V\to {\mathbb C}$, $b:V\times V\to{\mathbb C}$ be sesquilinear forms\footnote{One can often think of sesquilinear forms in terms of physical quantities. In fact, if $a,b$ are symmetric, $2a(u,u),2b(u,u)$ merely represent the energy functionals associated with the elastic and damping operators $A,B$ that appear in~\eqref{damped}, respectively.}.

More precisely, we recall that the~\emph{operator associated with $a$} is by definition given by
\begin{equation*}
\begin{array}{rcl}
D(A)&:=&\left\{f \in V:\exists h\in H \hbox{ s.t. } a(f,g)=(h\mid g)_H\; \forall g\in V\right\},\\
Af &:=&h,
\end{array}
\end{equation*}
and likewise for the operator associated with $b$.

\bigskip
The following perturbation lemma seems to be of independent interest. It is the form equivalent of a well-known perturbation result for operators due to Desch--Schappacher. In the following we denote by $H_\alpha$ any interpolation space between $V$ and $H$ that verifies the interpolation inequality
\begin{equation}\label{interpeq}
\|f\|_{H_\alpha}\leq M_\alpha\|f\|_V^\alpha \|f\|_H^{1-\alpha},\qquad f\in V.
\end{equation}

\begin{lem}\label{perturbd}
Let $a:V\times V\to\mathbb C$ be a sesquilinear mapping. Let $\alpha\in [0,1)$ such that $a_1:V\times H_\alpha\to \mathbb C$ and $a_2:H_\alpha\times V\to\mathbb C$ be continuous sesquilinear mappings. Then $a$ is $H$-elliptic if and only if $a+a_1+a_2:V\times V\to\mathbb C$ is $H$-elliptic.
\end{lem}

\begin{proof}
Let $a$ be $H$-elliptic and let
\begin{equation*}\label{contin}
|a_1(f,g)|\leq M \|f\|_V \|g\|_{H_\alpha}\qquad\hbox{and}\qquad |a_2(g,f)|\leq M \|g\|_V \|f\|_{H_\alpha},
\end{equation*}
for some constant $M>0$ and for all $f\in V$, $g\in H_\alpha$, so that by~\eqref{interpeq} we can estimate both 
$|a_1(f,f)|$ and $|a_2(f,f)|$ by $M M_\alpha \|f\|_V^{1+\alpha} \|f\|_H^{1-\alpha}.$

By Young's inequality one has for all $\alpha\in[0,1)$ and all $x,y>0$ that
$$xy\leq \frac{1+\alpha}{2}x^\frac{2}{1+\alpha} +\frac{1-\alpha}{2} y^\frac{2}{1-\alpha}.$$
Thus, for all $\epsilon>0$ letting $x=(\sqrt{\epsilon} \|f\|_V)^{1+\alpha}$ and $y=(\frac{1}{\sqrt{\epsilon}} \|f\|_H)^{1-\alpha}$ one obtains
$$\|f\|_V^{1+\alpha} \|f\|_H^{1-\alpha}\leq \frac{1+\alpha}{2}\epsilon \| f\|_V^2 +\frac{1-\alpha}{2\epsilon} \|f\|_H^2,\qquad f\in V.$$
Accordingly, for all $\epsilon>0$ there exists $M(\epsilon)>0$ such that
$$-\epsilon \|f\|^2_V +M(\epsilon)\|f\|^2_H\leq a_1(f,f)+a_2(f,f),\qquad f\in V.$$
By assumption $a$ is $H$-elliptic, i.e., ${\rm Re}a(f,f)\geq \alpha \| f\|^2_V -\omega\| f\|^2_H$ for some $\alpha>0$ and $\omega\in\mathbb R$. Thus, that for $\epsilon=\alpha/2$
\begin{eqnarray*}
{\rm Re}(a+a_1+a_2)(f,f)&=&{\rm Re}a(f,f)+a_1(f,f)+a_2(f,f)\\
&\geq& \alpha \| f\|^2_V -\omega\| f\|^2_H-\epsilon \|f\|^2_V -M(\epsilon)\|f\|^2_H\\
&\geq& \frac{\alpha}{2} \| f\|^2_V -(\omega+M(\epsilon))\| f\|^2_H,
\end{eqnarray*}
for all $f\in V$. This completes the proof.
\end{proof}

With the aim of discussing the abstract damped wave equation~\eqref{damped} we introduce ${\mathbf V}:=V\times V$ as well as the candidate energy space ${\mathbf H}:=V\times H$. Observe that ${\mathbf V}$ is continuously and densely imbedded into ${\mathbf H}$ and that both ${\mathbf V}$ and ${\mathbf H}$ have a canonical Hilbert space structure.
Define 
\begin{equation}\label{defn}
{\ea}({\eu},{\ev}):= -(u_2 \mid v_1)_V + a(u_1,v_2)+ b(u_2,v_2),
\end{equation}
where we have considered
$${\eu}=(u_1, u_2)^\top,{\ev}=(v_1, v_2)^\top \in {\mathbf V},$$
i.e., $\ea$ is a sesquilinear form with domain ${\mathbf V}$. Observe that $\ea$ is in general not symmetric.

\begin{lem}\label{closed} 
The following assertions hold.
\begin{enumerate}[1)]
\item The form ${\ea}$ is continuous with respect to ${\mathbf V}$ if and only if $a,b$ are continuous with respect to $V$.
\item The form ${\ea}$ is ${\mathbf H}$-elliptic if and only if $b$ is $H$-elliptic.
\item Let ${\rm Re}a(u,v)={\rm Re}(u\mid v)_V$ for all $u,v\in V$. If $b$ is accretive, then $\ea$ is accretive.
\item If  if $\ea$ is accretive, then $b$ is accretive.
\end{enumerate}
\end{lem}

Observe that, as a direct consequence of the sesquilinearity of $a$, ${\rm Re}a(u,v)={\rm Re}(u\mid v)_V$ for all $u,v\in V$ if and only if ${\rm Re}a(u,v)\geq {\rm Re}(u\mid v)_V$ for all $u,v\in V$.

\begin{proof}
% 2) Observe that the form norm of $\ea$ is given by
% \begin{eqnarray*}
% \Vert {\eu}\Vert^2_{\ea} &=& {\rm Re}{\ea}({\eu},{\eu})+ \Vert {\eu}\Vert^2_{\mathbf H}\\
% &=&\Vert u_1\Vert^2_V + \Vert u_2\Vert^2_b+{\rm Re} a(u_1,u_2)- {\rm Re}(u_2\mid u_1)_V,\qquad {\eu}=(u_1,u_2)^\top \in{\mathbf V}.
% \end{eqnarray*}
% Thus, by~\eqref{hyper} we obtain that 
% $$\Vert {\eu}\Vert^2_{\ea} \geq \Vert u_1\Vert^2_V + \Vert u_2\Vert^2_b,\qquad {\eu}=(u_1,u_2)^\top\in{\mathbf V}.$$
% If $\ea$ is closed, i.e., if ${\mathbf V}$ is closed in ${\mathbf H}$ with respect to $\Vert\cdot\Vert_{\ea}$, then it is clear that the same holds for $b$. 
% 
% 3) Conversely, let $b$ be closed and $a$ be continuous, i.e., 
% $$| a(u,v)| \leq M_a\Vert u\Vert \Vert v\Vert,\qquad u,v\in V,$$ 
% for some constant $M_a>0$. Then, by the Cauchy--Schartz inequality we obtain that 
% $${\rm Re} a(u,v)- {\rm Re}(v\mid u)_V\leq (M_a+1)\Vert u\Vert_V\Vert v\Vert_V,\qquad u,v\in V.$$
% We conclude that
% $$\Vert {\eu}\Vert^2_{\ea} \leq M_1(\Vert u_1\Vert^2_V + \Vert u_2\Vert^2_b),\qquad {\eu}=(u_1,u_2)^\top \in {\mathbf V},$$
% for some constant $M_1>0$. Thus, c) is proved.
% 
1) Let $\ea$ be continuous. Then for some constant $M_{\ea}>0$ and all $u,v\in V$ one has
$$| b(u,v)|=| \ea({\eu},{\ev})| \leq M_{\ea}\Vert {\eu}\Vert_{\mathbf V}\Vert {\ev}\Vert_{\mathbf V}=M_{\ea} \Vert u\Vert_V\Vert v\Vert_V,$$
where we have set ${\eu}:=(0,u)^\top$ and ${\ev}:=(0,v)^\top$. Similarly, setting ${\eu}:=(u,0)^\top$ and ${\ev}:=(0,v)^\top$ we obtain that
$$| a(u,v)|=| {\ea}({\eu},{\ev})| \leq M_{\ea}\Vert {\eu}\Vert_{\mathbf V}\Vert {\ev}\Vert_{\mathbf V}=M_{\ea} \Vert u\Vert_V\Vert v\Vert_V.$$

Let now $a,b$ be continuous, i.e., assume that for some $M_a,M_b\geq 0$ there holds
$$| a(u,v)| \leq M_a\Vert {u}\Vert_V \Vert v\Vert_V ,\qquad u,v\in V,$$
as well as
$$| b(u,v)| \leq M_b\Vert {u}\Vert_V \Vert v\Vert_V,\qquad u\in V.$$
A tedious computation then shows that
\begin{eqnarray*}
| {\ea}({\eu},{\ev})|^2 &\leq &
\Vert u_2\Vert_V^2 \Vert v_1\Vert_V^2+ M_a^2  \Vert u_1\Vert_V^2 \Vert v_2\Vert_V^2 + M_b^2\Vert u_2\Vert_V^2 \Vert v_2\Vert_V^2\\
&&\quad + 2M_a  \Vert u_1\Vert_V \Vert u_2\Vert_V \Vert v_1\Vert_V \Vert v_2\Vert_V 
+ 2  M_b\Vert u_2\Vert^2_V \Vert v_1\Vert_V \Vert v_2\Vert_V\\
&&\quad + 2M_a M_b   \Vert u_1\Vert_V\Vert u_2\Vert_V\Vert v_2\Vert_V^2\\
% &\leq & M_a^2  \Vert u_1\Vert_V^2 \Vert v_2\Vert_V^2 + \Vert u_2\Vert_V^2 \Vert v_1\Vert_V^2+	M_b^2\Vert u_2\Vert_V^2 \Vert v_2\Vert_V^2\\
% &&\quad + \frac{M_a  }{2}(\Vert u_1\Vert_V^2 + \Vert u_2\Vert_V^2)( \Vert v_1\Vert_V^2 +\Vert v_2\Vert_V^2) \\
% &&\quad 	+ M_a M_b  (\Vert u_1\Vert_V^2 + \Vert u_2\Vert_V^2)\Vert v_2\Vert_V^2\\
% &&\quad + M_b \Vert u_2\Vert_V^2 (\Vert v_1\Vert_V^2 + \Vert v_2\Vert_V^2)\\
% &\leq & \frac{M_a  }{2}\Vert u_1\Vert_V^2 \Vert v_1\Vert_V^2 +
%  (M^2_a   +\frac{M_a  }{2}+ M_a M_b  )\Vert u_1\Vert_V^2 \Vert v_2\Vert_V^2\\
% &&\quad  +(1+\frac{M_a  }{2}+ M_a M_b  )\Vert u_2\Vert_V^2 \Vert v_1\Vert_V^2 \\
% &&\quad + (M_b^2+\frac{M_a  }{2}+ M_a M_b   + M_b)\Vert u_2\Vert_V^2 \Vert v_2\Vert_V^2\\
&\leq & M_\ea^2(\Vert u_1\Vert_V^2 + \Vert u_2\Vert_V^2)(\Vert v_1\Vert_V^2 + \Vert v_2\Vert_V^2),
\end{eqnarray*}
i.e., $| {\ea}({\eu},{\ev})|\leq M_\ea \Vert \eu\Vert_{\mathbf V}\Vert {\ev}\Vert_{\mathbf V}$, where
\begin{equation}\label{angle}M_\ea^2:=\frac{M_a  }{2}+ M_a M_b  +\max\left\{M^2_a  ,1,M^2_b\right\}.
\end{equation}

2) To begin with, consider the form $\ea_0:{\mathbf V}\times {\mathbf V}\to{\mathbb C}$ defined by 
$$\ea_0(\eu,\ev):=b(u_2,v_2).$$ 
A direct computation shows that $\ea_0$ is $\mathbf H$-elliptic if and only if $b$ is $H$-elliptic.  Similarly, define the continuous sesquilinear mappings $\ea_1:{\mathbf H}\times {\mathbf V}\to{\mathbb C}$ and $\ea_2:{\mathbf V}\times {\mathbf H}\to{\mathbb C}$ by
$$\ea_1(\eu,\ev):=-(u_2\mid v_1)_V\qquad\hbox{and}\qquad \ea_2(\eu,\ev):=a(u_1,v_2).$$
By Lemma~\ref{perturbd} we conclude that $\ea=\ea_0+\ea_1+\ea_2$ is $\mathbf H$-elliptic if and only if $\ea_0$ is $\mathbf H$-elliptic if and only if $b$ is $H$-elliptic.

3) If $b$ is accretive and ${\rm Re}a(u,v)={\rm Re}(u\mid v)_V$ for all $u,v\in V$, then
$${\rm Re}{\ea}({\eu},{\eu})= {\rm Re}b(u_2,u_2)\geq 0,\qquad {\eu}= (u_1,u_2)^\top \in{\mathbf V},$$
i.e., $\ea$ is accretive. 

4) Conversely, if $\ea$ is accretive, we obtain that for all $u\in V$
$${\rm Re} b(u,u)={\rm Re}{\ea}({\eu},{\eu})\geq 0,$$
 where we have set ${\eu}:=(0,u)$.
\end{proof}

% \begin{rem}
% If~\eqref{hyper} is replaced by the stronger assumption 
% \begin{equation}\label{hyper2}
% {\rm Re} a(u,v)= {\rm Re}(v\mid u)_V,\qquad u,v\in V.
% \end{equation}
% then we can obtain a sharp estimate on the continuity constant of $\ea$. More precisely, if $a$ is c
% \end{rem}

By~\cite[Prop.~1.51 and Thm.~1.52]{Ou05} we can now state the following. 

\begin{thm}\label{gener}
Let $a,b$ be continuous. Let further $b$ be $H$-elliptic. Then the operator associated with $\ea$ is closed. It generates a $C_0$-semigroup $(e^{-t\ea})_{t\geq 0}$ on ${\mathbf H}$ which is analytic of angle $\frac{\pi}{2}-\arctan{M}$, where $M_{\bf a}$ is defined as in~\eqref{angle}. The semigroup $(e^{-t\ea})_{t\geq 0}$ is contractive if $b$ is accretive and ${\rm Re}a(u,v)={\rm Re}(v\mid u)_V$ for all $u,v\in V$.
\end{thm}

We emphasize that in the above theorem we are assuming $a$ neither to be $H$-elliptic, nor to be (quasi)accretive. In other words, the operator $A$ associated with $a$ need not be closed or (quasi)dissipative. Thus, in the limiting case of $A$ bounded from $D(B)$ to $H$, where $B$ is the operator associated with $b$, Theorem~\ref{gener} extends the well-posedness results of~\cite{CT89,XL98,CS05}. %Our estimate on the angle of analyticity of $(e^{-t\ea})_{t\geq 0}$ should be compared with~\cite[Thm.~6.4.2]{XL98}.
In this sense, we say that the leading term in~\eqref{damped} is not the elastic, but rather the damping one.

\begin{rem}\label{peccato}
1) Let $V\not=\{0\}$. The form $\ea$ is self-adjoint if and only if $b$ is self-adjoint and  $a(\cdot,\cdot)=-(\cdot \mid \cdot)_V$.
Let in fact ${\eu}:=(u,0)^\top$ and ${\ev}:=(0,v)^\top$, with $u,v\in V$, $v\not=0\not=u$. Then, one has
$${\ea}({\eu},{\ev})=a(u,v)\qquad\hbox{and}\qquad  {\ea}({\ev},{\eu})=-(v\mid u)_V.$$
On the other hand, if ${\eu}:=(0,u)^\top$ and ${\ev}:=(0,v)^\top$, with $u,v\in V$, $v\not=0\not=u$, then
$${\ea}({\eu},{\ev})=b(u,v)\qquad\hbox{and}\qquad  {\ea}({\ev},{\eu})=b(v,u).$$
To prove the converse implication, it suffices to observe that  if $b$ is self-adjoint and $a(\cdot,\cdot)=-(\cdot \mid \cdot)_V$ , then
$$\ea(\eu,\ev)=b(u_2,v_2)=\overline{b(v_2,u_2)}=\overline{\ea(\ev,\eu)}.$$

2) The form $\ea$ is not coercive, unless $V=\{0\}$. Let in fact ${\eu}:=(u,0)^\top$, with $0\not= u\in V$. Then one has
$${\rm Re}{\ea}({\eu},{\eu})=0<\Vert {\eu}\Vert^2_{\mathbf V}.$$
This shows that there exists no $\epsilon>0$ such that the estimate $\| e^{-t\ea}\|\leq e^{-\epsilon t}$ holds for all $t\geq 0$. This should be compared with the exponential stability result in~\cite[Thm.~1.1]{CT89}.

% 2) Assume that for all $u\in V$ one has 
% \begin{itemize}
% \item ${\rm Re}u\in V$,
% \item $a({\rm Re}u,{\rm Im}u),b({\rm Re}u,{\rm Im}u)\in {\mathbb R}$, and
% \item $a({\rm Re} u,{\rm Re}u)+\Vert {\rm Re}u\Vert^2_V \geq b({\rm Re}u,{\rm Re}u)$.
% \end{itemize}
% Then the form $\ea$ is not sectorial (in the sense of~\cite[Def.~1.7]{Ou05}) unless $V=\{0\}$. In fact, let $\gamma$ be an arbitrary positive constant. Let $w\in V$ such that $u_1:={\rm Re}w\not =0$, and set ${\eu}:=(\gamma u_1,iu_1)^\top$. Then
% $${\ea}({\eu},{\eu}):= -i\gamma \Vert u_1 \Vert^2_V - i \gamma a(u_1,u_1)+ b(u_1,u_1).$$
% Thus, one sees that
% $$| {\rm Im}{\ea}({\eu},{\eu})| = \gamma\Vert u_1 \Vert^2_V + \gamma a(u_1,u_1) > \gamma b(u_1,u_1)=\gamma {\rm Re}{\ea}({\eu},{\eu}).$$
% Therefore, the operator $\mathbf A$ is sectorial only upon translation of $\ea$.
% 

3) In the relevant case of ${\rm dim}\;{\mathbf V}=\infty$ the imbedding of ${\mathbf V}$ in ${\mathbf H}$ is not compact. Thus if Theorem~\ref{gener} applies, then $(e^{-t\ea})_{t\geq 0}$ is not compact.

% 4) We recall that the numerical range $W(c)$ of a sesquilinear, continuous, $H$-elliptic form $c$ is defined as the set
% $$W(c):=\{c(u,u)\in{\mathbb C}: \| u\|_H=1\}.$$
%  Then, if $b$ is accretive, a direct computation shows that 
%  $$W(b)\subset W(\ea)\subset W(b)+\{a(u,v)-(v\mid u)_V: \| u\|^2_V + \| v\|^2_H=1\}.$$
% It is known that several properties of the semigroup $(e^{-t\ea})_{t\geq 0}$ can be deduced by the location of  $W(\ea)$, cf.~\cite{Cr04}. This will be exploited in Section~3.
% However, if $(u,v)^\top \in V\times H$ with $\| u\|^2_V + \| v\|^2_H=1$, then by the continuity of $a$ one sees that $|a(u,v)-(v\mid u)_V|\leq 2\| u\|_V \| v\|_V\leq \frac{1}{2}(\| u\|^2_V + \| v\|^2_V)= \frac{1}{2}(1 -\| v\|^2_H + \| v\|^2_V)$. Thus, if $K>0$ is the embedding costant such that $\| u\|_H\leq K \| u\|_V$ for all $u\in V$, then $|a(u,v)-(v\mid u)_V|\leq \frac{1}{2}(1 + (1+K^2) \| v\|^2_V)$.
% 
% This estimate will be improved in the proof of Theorem~\ref{crocrit} below.
4) An advantage of dealing with sesquilinear forms instead of operators is the flexibility of this theory. Let us briefly discuss the case of time-dependent damped wave equations. Consider families $(a_{t})_{t\in [0,T]}$ and $(b_t)_{t\in [0,T]}$ of sesquilinear forms with joint (time-independent) dense domain $V$. Assume them to be equicontinuous. Let furthermore the mappings $t\mapsto a_t(u,v)$ and $t\mapsto b_t(u,v)$ be measurable for all $u,v\in V$. If finally $(b_t)_{t\geq 0}$ is equi-$H$-elliptic, i.e.,
$${\rm Re}b_t(u,u)+\omega \Vert u\Vert^2_H \geq \alpha \Vert u\Vert^2_V,\qquad u\in V,\; t\geq 0,$$
for some $\omega\in\mathbb R,\alpha>0$, then it is easy to see that the family of sesquilinear forms $({\ea}_t)_{t\geq 0}$ defined by
$${\ea}_t({\eu},{\ev}):= -(u_2 \mid v_1)_V + a_t(u_1, v_2)+ b_t(u_2,v_2),$$
fits the framework of Lions' theory of time-dependent forms, cf.~\cite{LM72}, and we conclude that the nonautonomous abstract Cauchy problem associated with $({\ea}_t)_{t\geq 0}$ is well-posed in a suitably weak sense.
\end{rem}

In order to interpret Theorem~\ref{gener} as a well-posedness result for~\eqref{damped}, we still have to determine the  operator $(\mathbf A,D({\mathbf A}))$ associated with  $\ea$, which  by definition is
\begin{eqnarray*}
D({\mathbf A})&:=&\{\eu\in {\mathbf V}: \exists {\ez}\in{\mathbf H}\hbox{ s.t. }\ea(\eu,\ev)=(\ez\mid \ev)_{\mathbf V} \hbox{ for all }\ev\in{\mathbf V}\},\\
{\mathbf A}\eu&:=&\ez.
\end{eqnarray*}

In fact, the expression ``$Au+B\dot{u}$" in~\eqref{damped} is in general purely formal, as the solution $u$ to~\eqref{damped} need not satisfy $u\in C({\mathbb R}_+,D(A))\cap C^1({\mathbb R}_+,D(B))$.
However, in our framework a direct computation shows that the following holds.

\begin{prop}\label{identAB}
The operator $\mathbf A$ on ${\mathbf H}$ associated with the form $\ea$ is given by
\begin{eqnarray*}\label{AB}
D({\mathbf A})&=&\{\eu\in{\mathbf V}: \exists w\in H\hbox{ s.t. } a(u_1,v)+b(u_2,v)=(w\mid v)_H \hbox{ for all }v\in V\},\\
{\mathbf A}\eu&=&(u_2,w)^\top.
\end{eqnarray*}
\end{prop}

% \begin{proof}
% By definition, the operator $\mathbf B$ associated with the form $\ea$ is given by
% \begin{eqnarray*}
% D({\mathbf B})&:=&\{{\eu}\in{\mathbf V}:\exists {\ez}\in {\mathbf H} \hbox{ s.t. } {\ea}({\eu},{\ev})=({\ez}\mid{\ev})_{\mathbf H}\; \forall{\ev}\in{\mathbf V}\}\\
% {\mathbf B}{\eu}&:=&-{\ez}.
% \end{eqnarray*}
% We want to show that ${\mathbf A}={\mathbf B}$. In order to prove that ${\mathbf A}\subset {\mathbf B}$, take ${\eu}=(u_1,u_2)^\top\in D({\mathbf A})$ and let ${\ez}:=(-u_2,-A(u_1+\rho u_2))^\top\in {\mathbf H}$. Then for all ${\ev}=(v_1,v_2)^\top\in{\mathbf V}$ one has
% $${\ea}({\eu},{\ev})= -(u_2 \mid v_1)_V + a(u_1+\rho u_2\mid v_2)=-(u_2 \mid v_1)_V - (A(u_1+\rho u_2)\mid v_2)_H=({\ez}\mid{\ev})_{\mathbf H},$$
% since $A$ is the operator associated with $a$.
% 
% Conversely, let ${\eu}=(u_1,u_2)^\top	\in D({\mathbf B})$ and ${\ez}=(z_1,z_2)^\top\in{\mathbf H}$ such that for all ${\ev}=(v_1,v_2)^\top\in{\mathbf H}$ one has
% $$-(u_2 \mid v_1)_V + (u_1+\rho u_2\mid v_2)={\ea}({\eu},{\ev})=({\ez},{\ev})_{\mathbf H}=(z_1\mid v_1)_V+ (z_2\mid v_2)_H$$
% Thus, $z_1=-u_2$ and by definition of operator associated with the form $a$ we deduce that $u_1+\rho u_2\in D(A)$ and moreover $z_2=-A(u_1+\rho u_2)$. This concludes the proof.
% \end{proof}

In the remainder of this section we assume $V,H$ to be function spaces over a measure space $(X,\mu)$. The following is a direct consequence of the above proposition and should be compared with the results of~\cite{CT90}.

\begin{cor}\label{abr}
Let $\rho\in H$ such that $\rho u\in V$ and $a(u,v)=b(\rho u,v)$ for all $u,v\in V$. Then
\begin{eqnarray*}
D({\mathbf A})&=&\{\eu\in{\mathbf V}: \exists w\in H\hbox{ s.t. } b(\rho u_1+u_2,v)=(w\mid v)_H \hbox{ for all }v\in V\}\\
&=&\{\eu\in{\mathbf V}: \rho u_1+u_2\in D(B)\},\\
{\mathbf A}\eu&=&(u_2, B(\rho u_1+u_2))^\top.
\end{eqnarray*}
where $B$ denotes the operator associated with $b$. 
\end{cor}

While throughout the paper we consider complex Hilbert spaces, it is of interest for applications to ensure that solutions to~\eqref{damped} are in fact real whenever the initial data are real. %In other words, we want to characterize the property of reality of $(e^{-t\ea})_{t\geq 0}$.
% Our method based on the theory of sesquilinear forms is also efficient whenever we wish to characterize some qualitative properties of the semigroup associated with $\ea$ in terms of analogous properties of the semigroup associated with $b$. Throughout this section we assume that $H=L^2(\Omega)$ for some $\sigma$-finite measure space $(\Omega,\mu)$, and further impose the Assumptions~\ref{cosine}.
 In the following we denote the closed convex subsets ${V}_{\mathbb R}$ and ${H}_{\mathbb R}$ defined by the real-valued functions belonging to $V$ and $H$, respectively.

\begin{prop}\label{real}
Let $a,b$ be continuous and $b$ be $H$-elliptic. Assume further that ${\rm Re}u\in V$ and moreover $a({\rm Re}u,{\rm Im}u), ({\rm Re}u\mid {\rm Im}u)_V\in\mathbb R$ for all $u\in V$. Then $(e^{-t\ea})_{t\geq 0}$ is real (i.e., it leaves invariant $V_{\mathbb R}\times H_{\mathbb R}$) if and only if the semigroup associated with $b$ is real (i.e., it leaves invariant $H_{\mathbb R}$).
\end{prop}

\begin{proof}
Without loss of generality we can assume both $b$ and $\mathbf a$ to be accretive, since reality of a semigroup is invariant under rescaling. 
Let the semigroup associated with $b$ be real. Then by~\cite[Prop.~2.5]{Ou05} one has ${\rm Re}u\in V$ for all $u\in V$ and $b({\rm Re}u,{\rm Im}u)\in{\mathbb R}$.
Thus, for an arbitrary ${\eu}=(u_1,u_2)^\top\in{\mathbf V}$, one has ${\rm Re}{\eu}=({\rm Re}{\eu_1},{\rm Re}{\eu_2})^\top\in {\mathbf V}$ and moreover
$${\ea}({\rm Re}{\eu},{\rm Im}{\eu})=-({\rm Re}u_2\mid {\rm Im}u_1)_V+a({\rm Re}u_1, {\rm Im}u_2)+ b({\rm Re}u_2,{\rm Im}u_2) \in{\mathbb R}.$$
Since the projection $\mathbf P$ of ${\mathbf H}$ onto ${\mathbf H}_{\mathbb R}$ is given by
$${\mathbf P}{\eu}=({\rm Re}u_1,{\rm Re}u_2),\qquad {\eu}=(u_1,u_2)^\top\in{\mathbf H},$$
the claim follows by~\cite[Thm.~2.2]{Ou05}.
% (Here we have used the fact that $({\rm Re}u_2\mid {\rm Im}u_1)_V\in\mathbb R$, which can be explained in the following way: $(\cdot\mid\cdot)_V$ defines a bounded sesquilinear form on the Hilbert space $V$. The associated operator is $I_V$, thus the induced semigroup is $(e^t I_V)_{t\geq 0}$, which is real. The claim follows by~\cite[Prop.~2.5]{Ou05}, if we consider $w:={\rm Im}u_1 + i{\rm Re}u_2$, $\tilde{w}:={\rm Im}u_2 + i{\rm Re}u_1$, and compute $({\rm Re}w\mid {\rm Im}w)_V$, $a({\rm Re}\tilde{w}, {\rm Im}\tilde{w})_V$.)
% 
Conversely, let $(e^{-t\ea})_{t\geq 0}$ be real and let $u\in V$. Set  ${\eu}:=(0,u)^\top\in{\mathbf V}$. Then, ${\rm Re}{\bf u}=(0,{\rm Re}u)\in {\bf V}$ and	 $b({\rm Re}u,{\rm Im}u)={\ea}({\rm Re}{\eu},{\rm Im}{\eu})\in{\mathbb R}$. 
\end{proof}

% \begin{exa}\label{bilapl}
% For $\rho\in {\mathbb C}$ define the forms
% $$a(u,v):=\rho \int_{\mathbb R^n}  \Delta u \overline{ v}\quad \hbox{ and }\quad b(u,v):=\int_{\mathbb R^n} \Delta u \overline{\Delta v}.$$
% Then, it follows by Proposition~\ref{identAB} that the operator associated with $\ea$ is
% \begin{eqnarray*}
% D({\mathbf A})&=&H^2({\mathbb R}^n)\times H^4({\mathbb R}^n),\\
% {\mathbf A}\eu&=&\left(u_2,\rho\Delta u_1+\Delta^2 u_2\right)^\top,
% \end{eqnarray*}
% and by Theorem~\ref{gener} the initial value problem associated with
% $$\ddot{u}(t,x)+\rho\Delta u(t,x)+\Delta^2 u(t,x)=0,\qquad t\geq 0,\; x\in{\mathbb R}^n,$$
% is governed by an analytic semigroup of angle $\frac{\pi}{2}-{\rm arctan}\frac{5}{2}$ on $H^2({\mathbb R}^n)\times L^2({\mathbb R}^n)$. By Proposition~\ref{real}, the solution $u$ is real-valued whenever $u(0,\cdot),\dot{u}(0,\cdot)$ are real-valued.
% \end{exa}

\section{Interpolation spaces and nonlinear problems}\label{interpo}

In Theorem~\ref{gener} we have shown that if $a,b$ are continuous and $b$ is $H$-elliptic, the form $\ea$ is associated with an analytic semigroup on $\mathbf H$. We can sharpen this result under the additional assumption that for some constant $M_b>0$
\begin{equation}\label{hacr}
| {\rm Im} b(u,u)| \leq M_b \Vert u\Vert_H \Vert u\Vert_V,\qquad u\in V.
\end{equation}

\begin{thm}\label{crocrit}
If~\eqref{hacr} holds, then the operator $\mathbf A$ associated with $\ea$ generates a cosine operator function on ${\mathbf H}$. Moreover, the form domain ${\mathbf V}$ is isometric to the fractional power domain $D(\lambda+{\mathbf A})^\frac{1}{2}$, for $\lambda>0$ large enough.
\end{thm}

\begin{proof}
We first show that $|{\rm Im}\ea(\eu,\eu)|\leq M_\ea \| \eu\|_{\mathbf V}\| \eu\|_{\mathbf H}$ for some constant $M_\ea$ and all $\eu\in{\mathbf V}$.
Let to this aim $\eu=(u_1,u_2)^\top \in {\mathbf V}$. Since $| a(u,v)| \leq M_a \Vert u\Vert_V \Vert v\Vert_V$ for some $M_a>0$ and all $u,v\in V$, there holds
\begin{eqnarray*}
| {\rm Im} \ea(\eu,\eu)|^2  &\leq& (1+M^2_a)\Vert u_1\Vert^2_V \Vert u_2\Vert^2_V + M_b^2 \Vert u_2\Vert_H^2 \Vert u_2\Vert_V^2\\
&&\quad +2M_a \Vert u_1\Vert_V^2 \Vert u_2\Vert_V^2 +2M_b(1+M_a) \Vert u_1\Vert_V \Vert u_2\Vert_H \Vert u_2\Vert_V^2\\
&\leq & \left((1 + M_a)^2 \Vert u_1\Vert_V^2 + M_b^2 \Vert u_2\Vert_H^2\right) \Vert u_2\Vert_V^2\\
&&\quad +M_b(1+M_a) \left(\Vert u_1\Vert_V^2 + \Vert u_2\Vert_H^2\right) \Vert u_2\Vert_V^2\\
&\leq & (1+M_a+M_b)^2\left(\Vert u_1\Vert_V^2 + \Vert u_2\Vert_H^2\right) \Vert u_2\Vert_V^2\\
&\leq &(1+M_a+M_b)^2 \Vert \eu\Vert_{\mathbf H}^2 \Vert \eu\Vert_{\mathbf V}^2.
\end{eqnarray*}
This shows in particular that the numerical range of $\ea$ is contained in a parabola (see~\cite[p.~204]{Ha06}) and thus, applying a  result due to Crouzeix~\cite{Cr04}, we promptly obtain that $\mathbf A$ generates a cosine operator function on ${\mathbf H}$.

Moreover, by Haase's converse of Crouzeix's theorem (see~\cite[\S~5.6.6]{Ar04}) there exists an equivalent scalar product $((\cdot\mid \cdot))_{\mathbf H}$ on ${\mathbf H}$ and $\lambda>0$ such that the numerical range  of $\ea_\lambda:=\ea +\lambda((\cdot\mid\cdot))_{\mathbf H}$ lies in a parabola. Now it follows by a result due to McIntosh (see again~\cite[\S~5.6.6]{Ar04}) that $\mathbf A$ has the square root property. This concludes the proof.
\end{proof}

The following result should be compared with~\cite[Thm.~XVIII.5.1]{DL88}.

\begin{cor}\label{wellp}
Let $B=B_0+B_1$, where $B_0$ is a self-adjoint and strictly positive definite operator. Assume $A$ to be bounded from $D(B_0^\frac{1}{2})$ to $D(B_0^{-\frac{1}{2}})$ and $B_1$ to be bounded from $D(B_0^\frac{1}{2})$ to $H$. Then problem~\eqref{damped} is governed by an analytic semigroup of angle $\frac{\pi}{2}$ on $D(B_0^{\frac{1}{2}})\times H$. 

In particular,~\eqref{damped} admits a unique mild solution for all initial data $u_{10}\in D(B_0^{\frac{1}{2}})$ and $u_{20}\in H$.
If $A=\rho B$ for $\rho\in\mathbb C$, then~\eqref{damped} admits a unique classical solution for all $u_{10},u_{20}\in D(B_0^{\frac{1}{2}})$ such that $\rho u_{10}+u_{20}\in D(B)$.
\end{cor}

\begin{proof}
Let $b_0:D(B_0)\times D(B_0)\to{\mathbb C}$ the coercive, symmetric sesquilinear form associated with $B_0$. In particular, $B_0$ has the square root property (cf.~\cite[\S~5.5.1]{Ar04}) and therefore the form norm of $b_0$ is isomorphic to $D(B_0^\frac{1}{2})$. Since now for the sesquilinear form $b$ associated with $B$ holds
$$|{\rm Im} b(u,u)|=|{\rm Im} (B_0 u\mid u)_H+{\rm Im} (B_1u|u)_H|\leq \| B_1 u\|_H \|u\|_H\leq M\|u\|_{D(B_0)} \|u\|_H$$
for some constant $M>0$, one sees that~\eqref{hacr} is satisfied. After defining by $a$ the sesquilinear form associated with $A$, Theorem~\ref{crocrit} can be applied. Since every cosine operator function generator also generates an analytic semigroup of angle $\frac{\pi}{2}$ (see~\cite[Thm.~3.14.17]{ABHN01}), the claim holds.
\end{proof}

\begin{exa}\label{nonmult}
For an open bounded domain $\Omega\subset {\mathbb R}^n$ with  $C^2$-boundary $\partial \Omega$ consider the  complete second order problem
\begin{equation*}
\left\{\begin{array}{rcll}
\ddot{u}(t,x)&=&\nabla\cdot\Big(\alpha(x)\nabla u(t,x)+\beta(x) \nabla \dot{u}(t,x)\Big),\qquad &t\geq 0,\; x\in\Omega,\\
u(t,z)&=&\dot{u}(t,z)=0,&t\geq 0,\; z\in\partial \Omega,\\
u(0,x)&=&u_{10}(x), &x\in \Omega,\\
\dot{u}(0,x)&=&u_{20}(x), &x\in \Omega,
\end{array}
\right.
\end{equation*}
where $\alpha,\beta\in C^1(\overline\Omega)$ such that $0< \beta(x)$ for all $x\in\overline{\Omega}$.

Let $B=-\nabla\cdot (\beta\nabla)$ and $A=-\nabla\cdot (\alpha\nabla)$ on $H:=L^2(\Omega)$, both with domain $H^2(\Omega)\cap H^1_0(\Omega)$. Accordingly introduce the forms
$$b(f,g):=\int_\Omega \beta(x) \nabla f(x)\overline{\nabla g(x)}\qquad\hbox{and}\qquad a(f,g):=\int_\Omega \alpha(x) \nabla f(x)\overline{\nabla g(x)}.$$
 Then $D(B^\frac{1}{2})=H^1_0(\Omega)$ and by Corollary~\ref{wellp} and Corollary~\ref{abr} one concludes that the operator	
\begin{eqnarray*}
D({\mathbf A})&=&\left\{(u_1,u_2)^\top\in (H^1_0(\Omega))^2: \alpha \nabla u_1+\beta \nabla u_2\in H^1_0(\Omega)\right\},\\
{\mathbf A}\eu&=&\left(u_2,\nabla \left(\alpha \nabla u_1+\beta \nabla u_2\right)\right)^\top.
\end{eqnarray*}
generates  on $H^1_0(\Omega)\times L^2(\Omega)$ an analytic semigroup of angle $\frac{\pi}{2}$. This semigroup is contractive if $\alpha\equiv 1$ (and more generally also whenever $\alpha>0$, up to considering weighted phase space). It
yields the solutions to the above problem, which are real valued whenever $u_{10} \in H^1_0(\Omega)$ and $u_{20}\in L^2(\Omega)$ are real valued. 

The analytical well-posedness of the above problem has been shown in~\cite{CT89}  only in the case of $\alpha$ strictly positive, whereas we allow for $\alpha$ to be a complex-valued function.
% \begin{itemize}
% \item ${f}$ is continuous,
% \item $f$ is globally Lipschitz continuous in the last four variables, and moreover 
% \item $u_{10},u_{20}\in H^1_0(\Omega)$ with $\alpha \nabla u_{10}+ \beta \nabla u_{20}\in H^1_0(\Omega)$.
% \end{itemize}
\end{exa}

We can now exploit the technique developed in~\cite[Chapt.~7]{Lu95} for semilinear parabolic problems, which heavily relies on interpolation theory.  In order to avoid technicalities,  we consider in the remainder of this section the special case of $A=\rho B$ for some $\rho\in\mathbb C.$
This case is relevant in many concrete contexts, e.g., whenever investigating semilinear strongly damped equations like the Klein--Gordon one, see e.g.~\cite{HW75,Av87,ENS87,BP01}.
As an example of a possible application, we formulate the following, which is a direct consequence of~\cite[Thm.~7.1.3 and 7.1.10]{Lu95}. More refined results, also yielding global well-posedness, can be obtained by applying further tools from~\cite[\S~7.2]{Lu95}.

\begin{cor}\label{lunardi}
Let $B$ satisfy the assumptions of Corollary~\ref{wellp}. Assume $G:[0,T]\times D(B^\frac{1}{2})\times D(B^\frac{1}{2})\to H$ to be a continuous mapping that is locally H\"older continuous with respect to the first variable and locally Lipschitz continuous with respect to the second and third ones. Then for small initial data  $u_{10},u_{20}\in D(B^\frac{1}{2})$
\begin{equation*}
\left\{\begin{array}{l}
\ddot{u}(t)+B(\rho u+\dot{u})(t)=G(t,u_1(t),\dot{u_2}(t)),\qquad t\in[ 0,T],\\
u(0)=u_{10},\qquad \dot{u}(0)=u_{20},
\end{array}
\right.
\end{equation*}
has a unique classical solution, locally in time.
%If moreover $u_{20}\in D(B^\frac{1}{2})$ and $\rho u_{10}+ u_{20}\in D(B)$, then such a solution is strict and global.
 \end{cor}
% 
% \begin{cor}\label{lunardi2}
% Assume that ${\mathbf F}:[0,T]\times V\times V\to H$ is a continuous function that is globally Lipschitz continuous in the second and third variables. Then the semilinear abstract Cauchy problem
% \begin{equation*}
% \left\{\begin{array}{rcl}
% \dot{\eu}(t)&=&{\mathbf A}\eu(t)+{\mathbf F}(t,u_1(t),\dot{u_2}(t)),\qquad t\in[ 0,T],\\
% \eu(0)&=&(u_{10},u_{20})^\top.
% \end{array}
% \right.
% \end{equation*}
% has a classical solution for all $u_{10}\in V$ and $u_{20}\in H$. If moreover $u_{20}\in V$ and $(u_{10},u_{20})^\top\in D({\mathbf A})$, then such a solution is strict.
%  \end{cor}

Theorem~\ref{crocrit} also allows to apply the theory developed in~\cite{CL94} for quasilinear parabolic problems, where determining interpolation spaces is a crucial step, too. A prototypical result is the following, which can be compared with~\cite[Thm.~5.1]{CS05}.

\begin{cor}
Let $D$ be a subspace of $H$ with $D\hookrightarrow V$. Let the mapping
$$B:V\times V\to{\mathcal L}(\{(u,v)^\top \in V\times V: \rho u+v\in D\},H)$$
be well-defined and locally Lipschitz continuous. Let $u_{10},u_{20}\in V$ and assume the operator $B(u_{10},u_{20})$ %with domain $D$ 
to satisfy the assumptions of Corollary~\ref{wellp} with $D(B(u_{10},u_{20})^\frac{1}{2})=V$. Then for all $f\in L^2({\mathbb R}_+, H)$ and all $g\in {\rm Lip}({\mathbb R}_+\times V, H)$ there exists $\tau>0$ such that the problem
\begin{equation*}\label{dampedquasi}
\left\{\begin{array}{l}
\ddot{u}(t)+B(u(t),\dot{u}(t))(\rho u(t)+\dot{u}(t))=f(t)+g(t,u(t)),\qquad t\in (0,\tau),\\
u(0)=u_{10},\qquad \dot{u}(0)=u_{20},
\end{array}\right.
\end{equation*}
has a solution $u\in H^2((0,\tau),H)\cap H^1((0,\tau),V)$ with $\rho u+\dot{u}\in L^2((0,\tau),D)$.
% Let the mapping ${\mathbf A}:V\times V\ni (x,y)^\top \mapsto (I,{\mathbf A}(x,y))\in{\mathcal L}(\{(u,v)^\top \in V\times V: \rho u+v\in D(B)\},H)$, where
% $(I,{\mathbf A}(x,y))(u,v)=(v,B(x,y)(\rho u+v))$
\end{cor}

\section{Dynamic boundary conditions}\label{dynam}

We introduce a new Hilbert space, which we denote $\partial H$ because in the applications we have in mind this is often a space of boundary values of functions in $H$. We also consider a bounded linear operator $L:V\to \partial H$ with dense range (in $\partial H$) and dense kernel (in $V$) and the two Hilbert spaces
$${\mathcal V}:=\{{\mathfrak f}:=(f_1,f_2,f_3)^T\in V\times V\times \partial H:Lf_2=f_3\}\qquad\hbox{and}\qquad {\mathcal H}:=V\times H\times \partial H.$$
It follows from~\cite[Lemma~5.6]{MR06} that $\mathcal V$ is dense in $\mathcal H$.

In recent years it has become increasingly clear that the right functional setting in order to discuss equations with dynamic boundary conditions is obtained by enlarging the state space of the corresponding equation with non-dynamic boundary conditions, see e.g.~\cite{FGGR02,AMPR03}.

We define the sesquilinear form $\mathfrak a$ with dense domain $\mathcal V$ by
$${\mathfrak a}({\mathfrak u},{\mathfrak v}):=-(u_2|v_1)_V+a(u_1,v_2)+b(u_2,v_2),$$
i.e., $\mathfrak a$ acts formally on the first two coordinates of the vectors in its domain just like $\mathbf a$.
To begin with, we deduce a result analogous to Lemma~\ref{closed}.

\begin{lem}\label{closed2} 
The following assertions hold.
\begin{enumerate}[1)]
\item The form ${\mathfrak a}$ is continuous with respect to ${\mathcal V}$ if and only if $a,b$ are continuous with respect to $V$.
\item The form ${\ea}$ is ${\mathbf H}$-elliptic if the sesquilinear form 
\begin{equation}\label{formn}
\{(f_2,f_3)^T\in V\times \partial H:Lf_2=f_3\}\ni ((f_2,Lf_2),(g_2,Lg_2))\mapsto b(f_2,g_2)\in\mathbb C 
\end{equation}
is $H\times \partial H$-elliptic.
\item Let ${\rm Re}a(u,v)={\rm Re}(u\mid v)_V$ for all $u,v\in V$. If the form introduced in~\eqref{formn} is accretive, then $\mathfrak a$ is accretive.
\end{enumerate}
\end{lem}

We can finally prove a generation result in this setting, too.

\begin{thm}\label{gener2}
Let $a,b$ be continuous. Let further the form introduced in~\eqref{formn} be $H\times\partial H$-elliptic. Then the operator associated with $\mathfrak a$ is closed. It generates an analytic $C_0$-semigroup $(e^{-t\ea})_{t\geq 0}$ on ${\mathcal H}$. This is contractive if the form in~\eqref{formn} is accretive is accretive and ${\rm Re}a(u,v)={\rm Re}(v\mid u)_V$ for all $u,v\in V$. If additionally~\eqref{hacr} holds, then the operator associated with $\mathfrak a$ generates a cosine operator function on ${\mathcal H}$, too.
\end{thm}

Identifying the operator associated with $\mathfrak a$ is generally difficult. However, this can be accomplished in many concrete cases by a suitable application of the Gau{\ss}--Green formulae.

\begin{exa}
Consider again the forms $a,b$ (this time with domain $H^1(\Omega)$) and the operators $A,B$ introduced in Example~\ref{nonmult}, for the sake of simplicity with $\alpha\in\mathbb C$ and $\beta\equiv 1$. Moreover, let $L$ be the trace operator. By Theorem~\ref{gener2} one concludes that the operator $\mathcal A$ defined by
\begin{eqnarray*}
D({\mathcal A})&:=&\left\{{\mathfrak f}:=(\alpha f_1,f_2,Lf_2)^\top\in (H^1(\Omega))^2\times L^2(\partial\Omega): \Delta(\alpha f_1+ f_2)\in L^2(\Omega), \; \partial_\nu (\alpha f_1+f_2)\in L^2(\partial\Omega)\right\},\\
{\mathcal A}{\mathfrak f}&:=&\left(f_2,\Delta(\alpha f_1+f_2), \partial_\nu (\alpha f_1+f_2)\right)^\top.
\end{eqnarray*}
(where $\partial_\nu$ denotes the normal derivative) generates  on $H^1(\Omega)\times L^2(\Omega)\times L^2(\partial\Omega)$ a cosine operator function and hence a (contractive) analytic semigroup of angle $\frac{\pi}{2}$. A direct computation shows that this semigroups solves the initial-boundary value problem
% \begin{equation*}
% \left\{\begin{array}{rcll}
% \dot{u}(t,x)&=&v(t,x),\qquad &t\geq 0,\; x\in\Omega,\\
% \dot{v}(t,x)&=&\Delta(\alpha u(t,x)+ v(t,x)),\qquad &t\geq 0,\; x\in\Omega,\\
% \dot{v}(t,z)&=&\partial_\nu (\alpha u(t,z)+ v(t,z)),\qquad &t\geq 0,\; z\in\partial\Omega,\\
% u(0,x)&=&u_{10}(x), &x\in \Omega,\\
% v(0,x)&=&u_{20}(x), &x\in \Omega,\\
% v(0,z)&=&u_{30}(z), &z\in \partial\Omega,
% \end{array}
% \right.
% \end{equation*}
% i.e.,
\begin{equation*}
\left\{\begin{array}{rcll}
\ddot{u}(t,x)&=&\Delta(\alpha u(t,x)+ \dot{u}(t,x)),\qquad &t\geq 0,\; x\in\Omega,\\
\ddot{u}(t,z)&=&\partial_\nu (\alpha u(t,z)+ \dot{u}(t,z)),\qquad &t\geq 0,\; z\in\partial\Omega,\\
u(0,x)&=&u_{10}(x), &x\in \Omega,\\
\dot{u}(0,x)&=&u_{20}(x), &x\in \Omega,\\
\dot{u}(0,z)&=&u_{30}(z), &z\in \partial\Omega.
\end{array}
\right.
\end{equation*}
Admittedly, this well-posedness result could also have been obtained combining the results of~\cite{AMPR03} and Corollary~\ref{wellp}.
\end{exa}

\begin{rem}
For the sake of simplicity, we have avoided to consider an additional term describing further dynamic processes on the boundary space $H$. However, such a generalization can easily be achieved by standard perturbation theory.
\end{rem}

\end{document}